\documentclass[leqno,12pt]{article}
\usepackage{latexsym}
 \usepackage{graphicx}
\usepackage{amsmath}
\usepackage{amssymb}
\usepackage{mathtools}
\usepackage{cite}
\usepackage{color}
\usepackage{subfig}
\usepackage[margin=1.3in, top=1.3in, bottom=1.3in]{geometry}
\usepackage{algorithm}
\usepackage{algorithmic}
\usepackage{hyperref}
\usepackage{appendix}

\newcommand{\nc}{\newcommand}
\nc{\nt}{\newtheorem}
\nt{thm}{Theorem}[section]
\nt{cor}[thm]{Corollary}
\nt{prop}[thm]{Proposition}
\nt{lem}[thm]{Lemma}
\nt{defn}[thm]{Definition}
\nt{rem}[thm]{Remark}
\nt{exa}[thm]{Example}
\nt{ass}[thm]{Assumption}
\nt{alg}[thm]{Algorithm}
\nt{conj}[thm]{Conjecture}
\nt{claim}[thm]{Claim}
\nt{oracle}[thm]{Oracle}
\nc{\ip}[2]{\mbox{$\langle #1,#2 \rangle$}}
\nc{\pf}{\noindent{\bf Proof\ \ }}
\nc{\finpf}{\hfill{$\Box$}\linespace}
\nc{\linespace}{\vspace
{\baselineskip} \noindent}
\nc{\R}{{\mathbf R}}
\nc{\Z}{{\mathbf Z}}
\nc{\X}{{\mathbf X}}
\nc{\F}{{\mathcal F}}
\nc{\oR}{\overline{\R}}
\nc{\M}{\mathcal M}
\nc{\e}{\epsilon}

\nc{\Rn}{{\mathbf R}^n}
\nc{\inT}{\mbox{\rm int}\,}
\nc{\cl}{\mbox{\rm cl}\,}

\def\tto{\;{\lower 1pt \hbox{$\rightarrow$}}\kern -12pt
           \hbox{\raise 2.8pt \hbox{$\rightarrow$}}\;}
\newenvironment{myequation}{\setcounter{equation}{\value{thm}}
   \begin{equation}}{\addtocounter{thm}{1}\end{equation}}

\nc{\bmye}{\begin{myequation}}
\nc{\emye}{\end{myequation}}

\begin{document}
\title{
Convex optimization on CAT(0) cubical complexes
}
\author{Ariel Goodwin
\thanks{CAM, Cornell University, Ithaca, NY.
\texttt{arielgoodwin.github.io} 
}
\and
Adrian S. Lewis
\thanks{ORIE, Cornell University, Ithaca, NY.
\texttt{people.orie.cornell.edu/aslewis} 
\hspace{2cm} \mbox{~}
Research supported in part by National Science Foundation Grant DMS-2006990.}
\and
Genaro L\'opez-Acedo
\thanks{Department of Mathematical Analysis -- IMUS, University of Seville, 41012 Seville, Spain
\texttt{glopez@us.es} 
}
\and
Adriana Nicolae\thanks{Department of Mathematics, Babe\c{s}-Bolyai University, 400084, Cluj-Napoca, Romania \hfill \mbox{}
\texttt{anicolae@math.ubbcluj.ro} 
}
}
\date{\today}
\maketitle

\begin{abstract}
We consider geodesically convex optimization problems involving distances to a finite set of points $A$ in a CAT(0) cubical complex.  Examples include the minimum enclosing ball problem, the weighted mean and median problems, and the feasibility and projection problems for intersecting balls with centers in $A$.  We propose a decomposition approach relying on standard Euclidean cutting plane algorithms.  The cutting planes are readily derivable from efficient algorithms for computing geodesics in the complex.
\end{abstract}
\medskip

\medskip


\section{Introduction:  optimization in Hadamard space}
Metric spaces in which points can always be joined by geodesics --- isometric images of real intervals --- support a wide variety of interesting convex optimization problems \cite{zhang-sra,bacak-convex}.  Convex optimization is best understood in {\em Hadamard spaces\/}:  complete geodesic spaces satisfying the {\em CAT(0) inequality}, which requires that squared-distance functions to given points are all 1-strongly convex.  Manifolds comprise the most familiar examples, including Hilbert space, hyperbolic space, and the space of positive-definite matrices with its affine-invariant metric.

From a computational perspective, in the special case of manifolds, gradient-based methods are available \cite{boumal2022intromanifolds}.  On the other hand, even with no differentiable structure, some Hadamard spaces allow efficient computation of geodesics.  Particularly interesting are phylogenetic  tree spaces \cite{billera} and more general CAT(0) cubical complexes \cite{ardila}.  In these spaces, geodesics can be computed in polynomial time \cite{owen-provan,hayashi}.  Along with phylogenetic models, applications include reconfigurable systems in robotics \cite{abrams-ghrist,ghrist-peterson}.

Unfortunately, optimization algorithms in the setting of a general Hadamard space are scarce and slow.  Methods based on alternating projections are sometimes available \cite{sims-cyclic}, but other general-purpose algorithms typically rely on some predetermined sequence of step sizes.  The quintessential example is the problem of computing the mean of a finite set $A$ in the geodesic space 
$(\X,d)$, which is the unique minimizer of the function
\bmye \label{distance-squared}
m_A(x) ~=~ \sum_{a \in A} d^2(a,x) \qquad (x \in \X).
\emye
To compute the mean, the only known general methods iteratively update the $n$th iterate $x_n$ to a point on the geodesic between $x_n$ and some point $a_n \in A$, chosen randomly \cite{sturm} or via some deterministic strategy \cite{holbrook,bacak-medians,lim-palfia}:  a standard example is 
\[
x_{n+1} ~=~ \frac{n}{n+1} x_n ~+~ \frac{1}{n+1}a_n.
\]
To compute the mean of the set $\{-1,+1\} \subset \R$ by this method, for example, starting at the point $x_1 = \frac{1}{2}$ and using $a_n = (-1)^n$, results in the slowly converging sequence $x_n = \frac{1}{2n}(-1)^{n+1}$.  Computation in practice confirms this slow convergence:  see \cite{massart} for experiments in the manifold of positive definite matrices, and \cite{brown-owen} for experiments in phylogenetic tree space.

In this work, we focus on the particular case of geodesically convex optimization on a CAT(0) cubical complex.  This setting has several features suggesting faster algorithmic possibilities.  First, as we have noted, computing geodesics is tractable.  Secondly, the underlying space, by definition, decomposes into Euclidean cubes.  Lastly, restricted to each cube, geodesically convex optimization problems reduce to standard Euclidean convex optimization.  Using these ingredients, we suggest a new cutting plane approach.  Our central technique uses geodesics in CAT(0) cubical complexes to derive Euclidean subgradients of distance functions restricted to individual cells in the complex.  We illustrate the new algorithm on a simple computational example.

\section{CAT(0) cubical complexes}
We begin with some standard definitions \cite{bridson}.
Let $(\X,d)$ be a metric space. A {\em geodesic path} is a distance-preserving mapping $\gamma:[0,l] \subseteq \R \to \X$. The image $\gamma([0,l])$ is called a {\em geodesic segment}. We say that $\X$ is a {\em (uniquely) geodesic space} if every two points in $\X$ are joined by a (unique) geodesic segment.  A function $f \colon \X \to \R$ is {\em convex} when its composition with every geodesic path is convex.

Given three points $x, y, z$ in a geodesic space, a {\em geodesic triangle} $\Delta = \Delta(x,y,z)$ is the union of three geodesic segments (its sides) joining each pair of points. A {\it comparison triangle} for $\Delta$ is a triangle $\Delta(\bar x,\bar y,\bar z)$ in $\R^2$ that has side lengths equal to those of $\Delta$.  A geodesic metric space is called {\em CAT(0)} if in every geodesic triangle, distances between points on its sides are no larger than corresponding distances in a comparison triangle. CAT(0) spaces are uniquely geodesic.  A {\em Hadamard space} is a complete CAT(0) space.  For any nonempty closed convex subset $C$ of a Hadamard space $\X$, every point in $\X$ has a unique nearest point in $C$. 

Let $\gamma:[0,l]\to \X$ and $\eta:[0,r]\to \X$ be two nonconstant geodesics paths in a CAT(0) space issuing at the same point $x = \gamma(0) = \eta(0)$. For $t \in [0,l]$ and $s \in [0,r]$, denote $\gamma_t = \gamma(t)$ and $\eta_s = \eta(s)$, and let $\Delta(\bar{\gamma_t}, \bar{x}, \bar{\eta_s})$ be a comparison triangle for $\Delta(\gamma_t, x, \eta_s)$. Then the angle $\angle \bar{\gamma_t} \bar{x} \bar{\eta_s}$ is a nondecreasing function of both $t$ and $s$, and the {\em Alexandrov angle} between $\gamma$ and $\eta$ is defined by 
\[
\angle(\gamma,\eta) ~=~ \angle \gamma_l x \eta_r ~=~  \lim_{t,s \searrow 0} \angle \bar{\gamma_t} \bar{x} \bar{\eta_s}.
\]

A {\em polyhedral cell} $C$ is a geodesic metric space isometric to the convex hull $\hat C$ of finitely many points in a Euclidean space $\R^n$: we usually identify $C$ and $\hat C$.  By a {\em face} of $C$, we mean a nonempty set that is either $C$ itself, or is the intersection of $C$ with a hyperplane $H$ such that $C$ belongs to one of the closed half-spaces determined by $H$.  The {\em dimension} of a face is the dimension of the intersection of all affine subspaces containing it.  The $0$-dimensional faces of a cell are called its {\em vertices}, and the $1$-dimensional faces are called its {\em edges}. A {\em polyhedral complex} is a set of polyhedral cells of various dimensions such that the face of any cell is also a cell of the complex, and the intersection of any two cells is either empty or a face of both.  The complex is {\em finite} if it consists of finitely many cells, and is {\em cubical} if each $n$-dimensional cell is isometric to the unit cube $[0,1]^n$.

Given two points $x,y$ in a polyhedral complex $\X$, the distance $d(x,y)$ is the infimum of lengths of piecewise geodesic paths joining $x$ to $y$. A {\em piecewise geodesic path} from $x$ to $y$ is an ordered sets of points $x_0 = x, x_1, \ldots, x_k = y$ in $\X$ such that for every $i \in \{1, \ldots, k\}$, there exists a cell $C_i$ with $x_{i-1}, x_i \in C_i$; its length is $\sum_{i=1}^{k} |x_{i-1} - x_i|$. If $\X$ is connected and finite, then $(\X,d)$ is a complete geodesic metric space. 

In what follows, we will, in general, consider finite cubical complexes. As proved by Gromov \cite{gromov}, a cubical complex is CAT(0) if and only if it is simply connected and satisfies the following ``link'' condition at each vertex. Consider the edges containing that vertex.  Any set of $k$ distinct such edges, each pair of which is contained in a common 2-dimensional cell, must consist of edges all contained in a common $k$-dimensional cell.  For more discussion, see \cite{hayashi}. 

\begin{exa}[A simple CAT(0) cubical complex] \label{adriana}
{\rm
Consider the space
\[
\X ~=~ \{ x \in [-1,1]^2 : x_1 \le 0~ \mbox{or}~ x_2 \le 0 \},
\]
with the distance induced by the Euclidean metric:  in other words, the distance between points is the Euclidean length of the shortest path between them in $\X$.  This space is a finite CAT(0) cubical complex, consisting of the three $2$-dimensional cells
\[
P_1 = [-1,0] \times [0,1], \qquad P_2 = [-1,0] \times [-1,0], \qquad P_3 = [0,1] \times [0,-1].
\]

}
\end{exa}

More generally than the example above, any simply connected subcomplex of $\Z^2$, the integer lattice cubing of $\R^2$, is CAT(0), because no three edges containing a vertex can, pairwise, be edges of common squares, so the link condition holds.  On the other hand, consider the cubical complex formed from the cube $[0,1]^3$ in $\R^3$ by taking just the three faces containing zero.  This simply connected cubical complex is not CAT(0) because it fails the link condition:  the edges connecting zero with the three standard unit vectors are pairwise contained in common squares, but no cell contains all three.  

\section{Decomposition}
Before describing a decomposition approach to convex optimization on cubical complexes, we first illustrate by considering the mean of three points in the simple space of Example \ref{adriana}.

\begin{exa}[A simple mean calculation] \label{simple2}
{\rm
In the CAT(0) cubical complex $\X$ described in Example \ref{adriana}, consider the set $A$ consisting of the standard unit vectors $e_1$ and $e_2$, along with $-e_1$.  To compute the mean of $A$, we must solve the underlying optimization problem (\ref{distance-squared}).  To do so, we decompose $\X$ into the union of three cells $P_1$, $P_2$, and $P_3$, and solve the restricted problem over each cell in turn:
\[
\min_{x \in P_1} \{ (1+|x|)^2 + |x-e_2|^2 + |x+e_1|^2 \}
\]
\[
\min_{x \in P_2} \{ |x-e_1|^2 + |x-e_2|^2 + |x+e_1|^2 \}
\]
\[
\min_{x \in P_3} \{ |x-e_1|^2 + (1+|x|)^2 + |x+e_1|^2 \}.
\]
Worth noting is that the first and third problems are not smooth.  The optimal solution of both the second and third problems is the point $(0,0)$.  However, the optimal solution of the first problem is the strictly better point $(-\alpha,\alpha)$, where 
$\alpha = \frac{1}{6}(2-\sqrt{2})$.  This point is therefore the mean.
}
\end{exa}

This example illustrates a significant feature of mean calculations in cubical complexes or other Hadamard spaces that are not manifolds.  Even when we can compute geodesics efficiently, that tool alone does not immediately allow us to recognize whether or not a given point is the mean, let alone compute the mean.  For example, along the geodesics between the point $(0,0)$ and each point in the set $A$, the objective $m_A$ defined by equation (\ref{distance-squared}) is minimized at $(0,0)$.  However, $(0,0)$ is not the mean:  it does not minimize $m_A$ over the whole space $\X$. 

Returning to our general problem, we can minimize a convex function $f$ over a finite cubical complex $\X$ by minimizing $f$ over each of the finitely many cells comprising $\X$ separately.  To each of these subproblems we can apply a standard algorithm for Euclidean convex minimization.  
Rather than exhaustively optimizing over every cube, we can instead consider the following conceptual  method, inspired by a somewhat analogous conceptual approach sketched in \cite[Algorithm 4.4]{miller-owen-provan}.

\begin{alg}[Minimize convex $f$ on cubical complex $\X$] \label{algorithm}
{\rm
\begin{algorithmic}
\STATE
\STATE 	{\bf input:} initial point $x \in \X$
\STATE	{${\mathcal P} = \emptyset$}			\hfill \% {\em set of optimized cells}
\FOR{$\mbox{iteration} = 1,2,3,\ldots$ }
\STATE	{${\mathcal Q} = \{ \mbox{cells}~ P \not\in {\mathcal P} : x \in P \}$}
												\hfill \% {\em unoptimized cells containing $x$}
\IF{${\mathcal Q} = \emptyset$}
\RETURN	$x$										\hfill \% {\em $x$ optimal}
\ENDIF
\STATE	choose $P \in {\mathcal Q}$				
\STATE	choose $x_P$ minimizing $f$ over $P$	\hfill \% {\em solve new subproblem}
\IF{$f(x_P) < f(x)$}
\STATE	$x=x_P$									\hfill \% {\em best point so far}
\ENDIF
\STATE	${\mathcal P} = {\mathcal P} \cup \{\mbox{faces of}~P\}$	\hfill \% {\em update set of optimized cells}
\ENDFOR
\end{algorithmic}
}
\end{alg}

\noindent
While unnecessary formally, we would naturally always choose a {\em maximal} cell $P$, meaning that no strictly larger cell contains $x$.  

\begin{prop}[Termination]
For any finite cubical complex $\X$ and any convex function $f \colon \X \to \R$, Algorithm \ref{algorithm} terminates, returning a minimizer of $f$.
\end{prop}

\pf
The procedure terminates, since the space $\X$ is a finite complex.  Throughout the procedure, at the current iterate $x$, the value $f(x)$ never increases.  At termination, $x$ therefore minimizes the objective $f$ over every cell containing $x$.  Since $f$ is convex, $x$ therefore minimizes it over the whole space $\X$.
\finpf

As an example, consider the simple example above, on the space (\ref{adriana}).  In Algorithm \ref{algorithm}, as soon as we choose the cube $P_1$, the procedure terminates with the iterate $x$ equal to the mean, whether or not $P_2$ or $P_3$ have already been searched.

\subsection*{Means in metric trees}
It is illuminating to consider the behavior of Algorithm \ref{algorithm} for computing means for the simplest class of CAT(0) cubical complex:  the case when the space $\X$ is a {\em finite metric tree}.  In that case, $\X$ consists of the edges and vertices of a finite connected acyclic graph, where we identify edges with \mbox{1-dimensional} cells of unit length, intersecting at common 0-dimensional cells --- the vertices.  The Gromov link condition holds trivially.  

We consider the problem of computing the mean of a finite subset $A$ of finite metric tree.
We illustrate with the following example \cite[Example 1]{sticky}, which is a special case of the ``open books'' \cite{sticky} discussed in Appendix \ref{appendix}.

\begin{exa}[Stickiness]
{\rm
The {\em 3-spider} is the finite metric tree consisting of three copies of the interval $[0,1]$ joined at the shared origin $0$.  Denote the three copies $P^1,P^2,P^3$, and consider a set $A$ consisting of three points $a^i \in \inT P^i$ satisfying
\[
a^i < \sum_{j \ne i} a^j \qquad \mbox{for}~ i=1,2,3.
\]
(One particular example is $a^i = \frac{1}{2}$ for each $i$.)
A quick calculation shows that the mean is the shared origin $0$.  Unlike in the Euclidean case, the mean is insensitive to small changes in the points $a^i$:  it is {\em sticky} in the sense of \cite{sticky}.
}
\end{exa}

Consider a general finite metric tree $\X$, with vertex set $V$ and edge set $E$. Via a finite computation, we can exactly compute the mean of a finite subset $A$ using Algorithm~\ref{algorithm} to minimize the function $m_A$ in equation (\ref{distance-squared}).  The data of the problem, in addition to the graph $(V,E)$, consists, for each point $a \in A$, of an endpoint $v_a \in V$ of an edge $e_a \in E$ containing $a$, and the distance $\delta_a$ between $a$ and $v_a$ along the edge $e_a$.  

Problems of this kind have been widely studied in the operations research literature, in the context of facility location problems \cite{hansen}.  More typical than the mean problem in that context are the ``median'' or ``minimax'' problems, involving the sum or maximum of the distance functions rather than sum of their squares, although \cite[Section 3.4]{hansen} notes an efficient algorithm for the mean problem due to Goldman.  Goldman's 1972 algorithm \cite{goldman} for the ``1-center problem on a tree network'', in the terminology of \cite{tansel-francis-lowe}, involves a ``trichotomy'' at each iteration:  after checking an edge the algorithm stops or is confined to one of the two subtrees resulting from deleting that edge.  

Algorithm \ref{algorithm}, it transpires, has the same trichotomy property.  At the outset of each iteration, we have a current vertex $v \in V$, and a current set 
${\mathcal E} \subset E$ of already optimized edges.  For each point $a \in A$, we first find the unique sequence $P_a$ of edges joining $v$ to $v_a$:  its cardinality $|P_a|$ is the distance $d(v,v_a)$.   The distance from $a$ to $v$ is therefore given by
\[
d(a,v) ~=~
\left\{
\begin{array}{ll}
|P_a| + \delta_a & (e_a \not\in P_a) \\
|P_a| - \delta_a & (e_a \in P_a).
\end{array}
\right.
\] 

We next choose a vertex $v'$ neighboring the vertex $v$ and with corresponding edge $e=vv'$ outside the set ${\mathcal E}$, terminating if there is no such edge.  We then minimize $m_A$ over $e$.  Denote the set of those points $a \in A$ for which $e \in P_a$ by $\hat A$.  If we identify $e$ with the unit interval $[0,1]$, where $v$ corresponds to the point 0 and $v'$ corresponds to the point 1, then for any point $x$ on $e = [0,1]$, we have
\[
d(a,x) ~=~ 
\left\{
\begin{array}{ll} 
d(a,v) + x & (e \not\in P_a) \\
d(a,v) - x & (e \in P_a).
\end{array}
\right.
\]
We now find the unique point $x \in [0,1]$ minimizing the strictly convex function
\[
\sum_{a \in \hat A} (d(a,v) - x)^2 ~+~  \sum_{a \not\in \hat A} (d(a,v) + x)^2.
\]
The unconstrained minimizer is
\[
\bar x ~=~ \frac{1}{|A|} \Big( \sum_{a \in \hat A} d(a,v) -  \sum_{a \not\in \hat A} d(a,v) \Big).
\]
If $0< \bar x  < 1$, then the algorithm terminates:  the mean is the point on $e$ at a distance $\bar x$ from $v$.  Otherwise, we update ${\mathcal E}$ to include $e$, update the current iterate to $v'$ if 
$\bar x \ge 1$, and repeat.  

\subsection*{Solving the subproblems}
In the case of mean computations for metric trees, Algorithm \ref{algorithm} involves one-dimensional  subproblems with closed-form solutions.  In general, however, the subproblems are multivariate, requiring iterative techniques.  

The approach outlined in \cite[Algorithm 4.4]{miller-owen-provan} for computing means in the phylogenetic tree space of \cite{billera} relies on a smooth but nonconvex interior-penalty philosophy, the complexity of which is unclear.  By contrast, our approach via Algorithm \ref{algorithm} generates subproblems that, while nonsmooth, are convex.  At each iteration of Algorithm \ref{algorithm} we can identify the cell $P$ isometrically with a cube $[0,1]^n$, for some dimension $n$, and then apply a linearly convergent Euclidean cutting plane algorithm, efficient in theory and reliable in practice.  

One such cutting plane approach for general convex objectives $f$ would be to apply the randomized method of \cite{lee-sidford-vempala}, which needs just $\tilde O(n^2)$ evaluations of $f(x)$ to approximately solve the minimization problem over the cube $[0,1]^n$.  More precisely, if $\max_P f - \min_P f \le 1$, then with constant probability the excess $f(x) - \min_P f$ is reduced to $\epsilon > 0$ after no more than $O\big(n^2\log^{O(1)}(\frac{n\sqrt{n}}{\epsilon})\big)$ function evaluations.  In essence, the method is subgradient-based:  at each iteration, it approximates a subgradient, using $\tilde O(n)$ function evaluations.  

In this work, however, we are primarily interested in structured functions $f$ composed simply from distance functions:  a typical example is the mean objective $m_A$ in equation (\ref{distance-squared}).  Our central observation is that such objectives support methods based on explicit subgradients, an approach with three potential advantages.  First, we arrive at a deterministic rather than randomized algorithm.  Secondly, the complexity of the available algorithms is better:  the classical ellipsoid method still requires $\tilde O(n^2)$ function evaluations, but more recent cutting plane algorithms improve this to $\tilde O(n)$ \cite{vaidya}.  Lastly, we can experiment with algorithms known to be effective in practice, like proximal or level bundle methods \cite{preamb}.

\section{Subgradients of distance functions}
In a CAT(0) cubical complex $\X$, the geodesic between any two points, $a$ and $x$, is computable
in polynomial time \cite{ardila}.  Suppose that $x$ lies in a cell $P$, and consider the distance function to $a$, restricted to $P$.  We next show how to use the geodesic to calculate a subgradient of this function at the point $x$. 

We start with a simple tool.

\begin{lem}\label{cat(0)-triangle}
Let $\X$ be a CAT(0) space and consider three points $a, x, w \in \X$ with $a \ne x$ and $w \ne x$. Then
\[
\frac{d(a,x) - d(a,w)}{d(x,w)} ~\le~ \cos(\angle axw).
\]
\end{lem}

\pf
From the triangle inequality we know
\[
|d(a,x) - d(a,w)| ~\le~ d(x,w).
\]
We deduce
\[
d^2(a,x) - 2d(a,x)d(a,w) + d^2(a,w) ~\le~ d^2(x,w),
\]
and hence
\[
2d(a,x)\big( d(a,x) - d(a,w) \big) ~\le~ d^2(a,x) + d^2(x,w) - d^2(a,w).
\]
Consequently we have
\[\frac{d(a,x) - d(a,w)}{d(x,w)} ~\le~ \frac{d(a,x)^2 + d(x,w)^2 - d(a,w)^2}{2d(a,x)d(x,w)} ~\le~ \cos(\angle axw),
\]
where the second inequality follows from the law of cosines.
\finpf

\begin{thm} \label{exit}
In a CAT(0) cubical complex $\X$, consider two cells $P$ and $Q$ with a common face $F$, and points 
$x \in F$ and $y \in Q$.  Denote the nearest point to $y$ in $P$ by $z$.  Then $z$ is also the nearest point to $y$ in $F$. Furthermore, for any point $w \in P$, if $w \ne x$ and $x \ne z$, then the following angle inequality holds: 
\bmye \label{cos-ineq}
\cos(\angle yxw) ~\le~  \cos(\angle yxz) \cdot \cos(\angle zxw).
\emye 
\end{thm}

Before proving the theorem, it helps to keep in mind a simple example.

\begin{exa} \label{3d}
{\rm
Consider any cubical subcomplex $\X$ of $\Z^3$ containing the cells
\[
P ~=~ [0,1] \times [0,1] \times [0,1] \qquad \mbox{and} \qquad Q ~=~ [-1,0] \times [0,1] \times \{0\}.
\]
Their common face is the cell
\[
F ~=~ \{0\} \times [0,1] \times \{0\}.
\]
Consider the two points
\[
x ~=~ \Big(0,\frac{1}{3},0\Big) ~\in~ F\qquad \mbox{and} \qquad 
y ~=~ \Big(-1,\frac{2}{3},0\Big) ~\in~ Q.
\]
The closest point to $y$ in $P$, in either the Euclidean distance or the distance it induces in $\X$, is
\[
z ~=~  \Big(0,\frac{2}{3},0\Big) ~\in~ F.
\]
At $x$, the angle between the nontrivial geodesic $[x,y]$ and any nontrivial geodesic $[x,w]$ with 
$w \in P$ is minimized (either in $\X$ or in $\Rn$) when $w = z$. In fact, the stronger inequality \eqref{cos-ineq} holds.
}
\end{exa}

\noindent
{\bf Proof of Theorem \ref{exit}}.
The result only involves geodesics with endpoints in $P$ and $Q$. All such geodesics lie in the interval between $P$ and $Q$, in the terminology of \cite{ardila}, so without loss of generality we can suppose by \cite[Theorem 3.5]{ardila} that $\X$ is a subcomplex of $\Z^n$, the integer lattice cubing of $\Rn$.

Without loss of generality, we can suppose that $P$ and $Q$ are unit cubes in $\Z^n$ that are contained in $[-1,1]^n$ and both contain $0$.  Corresponding to any two partitions of the index set $\{1,2,3,\ldots,n\}$ into disjoint subsets,
\[
I_\le \cup I_= \cup I_\ge \qquad \mbox{and} \qquad I^\le \cup I^= \cup I^\ge,
\]
we can define such cubes in $[-1,1]^n$ using the relationship set $\Omega = \{\le,=,\ge\}$ by
\begin{eqnarray*}
u \in P & \Leftrightarrow &  \quad u_i \sim 0~ \mbox{for all}~ i \in I_\sim ~ \mbox{and}~ \sim\: \in \Omega \\
u \in Q & \Leftrightarrow &  \quad u_i \sim 0~ \mbox{for all}~ i \in I^\sim ~ \mbox{and}~ \sim\: \in \Omega.
\end{eqnarray*}
We lose no generality in assuming that $P$ and $Q$ are of this form.  To illustrate, Example \ref{3d} corresponds to the partitions 
\begin{eqnarray*}
I_\le = \emptyset \quad & I_= = \emptyset &  \quad I_\ge = \{1,2,3\} \\
I^\le = \{1\} \quad & I^= = \{3\}  & \quad I^\ge = \{2\}.
\end{eqnarray*}

The point $x$ lies in the common face $F$, which is the cube in $[-1,1]^n$ defined by
\[
u \in F \quad \Leftrightarrow \quad
u_i ~
\left\{
\begin{array}{ll}
\le 0 & \mbox{if}~ i \in I_\le \cap I^\le \\
\ge 0 & \mbox{if}~ i \in I_\ge \cap I^\ge \\
=   0 & \mbox{otherwise}.
\end{array}
\right.
\]
Since $y \in Q$ we deduce
\[
 y_i \sim 0~ \mbox{for all}~ i \in I^\sim ~ \mbox{and}~ \sim\: \in \Omega.
\]
It is easy to compute componentwise the {\em Euclidean} nearest point $z'$ to the point $y$ in the cube $P$:  since $y \in [-1,1]^n$, for each $i=1,2,3,\ldots,n$ we have
\begin{eqnarray*}
z'_i &=&
\left\{
\begin{array}{ll}
\min\{y_i,0\}	& (i \in I_\le)	\\
0				& (i \in I_=)	\\
\max\{y_i,0\}	& (i \in I_\ge)
\end{array}
\right. \\ \\
 &=&
\left\{
\begin{array}{ll}
y_i	& \mbox{if $i \in I_\le$ and $y_i \le 0$, or if $i \in I_\ge$ and $y_i \ge 0$}	\\
0				& \mbox{otherwise}.	
\end{array}
\right.
\end{eqnarray*}
Since $y \in Q$, we deduce
\[
z'_i
~=~
\left\{
\begin{array}{ll}
y_i	& \mbox{if $i \in \bar I$}	\\
0				& \mbox{otherwise},	
\end{array}
\right.
\]
for the index set 
\[
\bar I ~=~ (I_\le \cap I^\le) \cup (I_\ge \cap I^\ge).
\]
Clearly we have $z' \in F \subset Q$, and hence the line segment $[z',y]$ is contained in $Q$.  That line segment is therefore also a geodesic in $\X$.  Since distances between points in $\X$ are no less than their Euclidean distance, we deduce $z'=z$.

Since $[x,w] \subset P$ and $[x,y] \subset Q$, those line segments are also geodesics in $\X$.  The distance in $\X$ between pairs of points, one on each line segment, is never less than the Euclidean distance, so in $\X$, the angle between those geodesics is never less than their Euclidean angle. Moreover, since the points $y$, $x$, and $z$ all lie in the cube $Q$, and the points $z$, $x$, and $w$ all lie in the cube $P$, the angles $\angle yxz$ and $\angle zxw$ in $\X$ equal the corresponding Euclidean angles. Thus, it is enough to prove \eqref{cos-ineq} for Euclidean angles, and hence, in what follows, we consider all angles to be Euclidean.

Since $x \in F$, the formula above for $z$ shows $\angle yzx = \pi/2$.  Consequently we deduce $|x-y|^2 = |x-z|^2 + |z-y|^2$ and
\[\cos(\angle yxz) ~=~ \frac{|x-z|}{|x-y|}.\]
 Applying the law of cosines, we get
\[\cos(\angle yxw) ~=~ \frac{|x-y|^2 + |x-w|^2 - |w-y|^2}{2|x-y||x-w|}\]  
and
\[\cos(\angle zxw) ~=~ \frac{|x-z|^2 + |x-w|^2 - |w-z|^2}{2|x-z||x-w|}.\]
Thus, 
\[\cos(\angle yxw) - \cos(\angle yxz) \cos(\angle zxw) ~=~ \frac{|z-y|^2 + |w-z|^2 - |w-y|^2}{2|x-y||x-w|} ~\le~ 0,\]
where the last inequality follows because $z$ is the Euclidean nearest point to $y$ in $P$ and $w \in P$.
\finpf

Consider any distinct points $x$ and $a$ in a CAT(0) cubical complex $\X$.  The CAT(0) property ensures that the distance function to $a$, denoted by $d_a$ is a convex function on $\X$ \cite[Example 2.2.4]{bacak-convex}.  Using the algorithm of \cite{ardila}, we can compute the geodesic $[x,a]$:  it passes via a sequence of nontrivial geodesics
\[ 
[x,y_0]~,~[y_0,y_1]~,~[y_1,y_2]~,\ldots,~ [y_m,a]
\]
through some corresponding connected sequence of cells $Q_0,Q_1,Q_2,\ldots,Q_{m+1}$ (neither sequence being necessarily unique).  We refer to $[x,y_0]$ as an {\em initial segment} and to $Q_0$ as a corresponding {\em initial cell}.

Given a cube $P \subset \Rn$, a vector $v \in \Rn$ is a {\em subgradient} of a convex function 
$e \colon P \to \R$ at a point $x \in P$ if
\[
\ip{v}{w-x} \le e(w) - e(x) \qquad \mbox{for all}~ w \in P.
\]
The {\em subdifferential} $\partial e(x)$ is just the set of such subgradients. 

We can now derive our central tool, which allows us to compute subgradients of distance functions restricted to cells in cubical complexes. 

\begin{thm}[Subgradients of distance functions] \label{main}
Consider a cube $P \subset \Rn$ and a point $x \in P$.  Suppose that $P$ is a cell in a CAT(0) cubical complex $\X$, and consider any point $a \in \X$.  Let $d_a^P$ denote the restriction of the distance function $d_a$ to $P$.  If $a=x$ then  $0 \in \partial d_a^P(x)$.  Suppose, on the other hand, that $a \ne x$.  In the geodesic $[x,a]$, let $[x,y]$ be an initial segment corresponding to an initial cell $Q$.  
Denote by $F$ the common face of $P$ and $Q$.  In any ambient Euclidean space for $Q$, denote the nearest point to $y$ in $F$ by $z$.  If $z=x$, then $0 \in \partial d_a^P(x)$.  On the other hand, if $z \ne x$, then
\[
\frac{\cos(\angle yxz)}{|x-z|} (x-z) ~\in~ \partial d_a^P(x).
\]
\end{thm} 

\pf
The case $a=x$ is trivial, so assume henceforth $a \ne x$.  By definition, we have $y \ne x$.  

Consider first the case $z=x$.  In that case, by Theorem~\ref{exit}, $x$ is the nearest point to $y$ in $P$, and since $y \in [x,a]$, we deduce that $x$ is also the nearest point to $a$ in $P$ (see \cite[Proposition II.2.4]{bridson}), so $0 \in \partial d_a^P(x)$ follows.

On the other hand, consider the case $z \ne x$.  We must prove
\[
\frac{\cos(\angle yxz)}{|x-z|}\ip{x-z}{w-x} ~\le~ d^P_a(w) - d^P_a(x) 
\qquad \mbox{for all}~ w \in P.
\]
This inequality holds when $w=x$.  When $w \ne x$, we have
\[
\frac{d(a,x) - d(a,w)}{d(x,w)} ~\le~ \cos(\angle axw)  ~=~  \cos(\angle yxw)
~\le~  \cos(\angle yxz) \cdot \cos(\angle zxw),
\]
by Lemma \ref{cat(0)-triangle} and Theorem \ref{exit}.  The result now follows.
\finpf

To illustrate the result, we consider a simple example.

\begin{exa}[Calculating a subgradient of a distance function] \mbox{} \\
{\rm
Consider the CAT(0) cubical complex $\X$ the cells of which are the following squares in $\R^2$, along with their edges and vertices:
\begin{eqnarray*}
P = [0,1] \times [0,1] & & Q = [0,1] \times [-1,0] \\
S = [0,1] \times [-2,-1] & & T = [-1,0] \times [-2,-1].
\end{eqnarray*}
For the point $a = (-\frac{1}{2},-2) \in \X$, the distance function $d_a \colon \X \to \R$ is given by
\[
d_a(w) ~=~ 
\left\{
\begin{array}{ll}
\frac{1}{2}\sqrt{5} + \sqrt{w_1^2 + (w_2 + 1)^2}	& (w_1 \ge 0~\mbox{and}~ w_2 \ge 2w_1 - 1) \\ \\
\sqrt{(w_1 + \frac{1}{2})^2 + (w_2 + 2)^2}			& (\mbox{otherwise}).
\end{array}
\right.
\]
This function must be convex on $\X$, so in particular it is convex on $P \cup Q$.  Differentiating shows that, for points $w$ in the interior of the set $\X$, regarded as a subset of 
$\R^2$, we have
\[
\nabla d_a(w) ~=~ 
\left\{
\begin{array}{ll}
\frac{1}{\sqrt{w_1^2 + (w_2 + 1)^2}} (w_1,w_2 + 1)	& (w_1 > 0~\mbox{and}~ w_2 > 2w_1 - 1) \\ \\
\frac{1}{\sqrt{(w_1 + \frac{1}{2})^2 + (w_2 + 2)^2}} (w_1 + \frac{1}{2}, w_2 + 2) & (w_1 < 0~ \mbox{or} ~ w_2 < 2w_1 - 1).
\end{array}
\right.
\]
The two cases in this formula describe a set of full measure around the point $x = (\frac{1}{2},0)$, and as $w$ approaches $x$ within this set, the gradient has a unique limit:  $\frac{1}{\sqrt{5}}(1,2)$.  Standard properties of Euclidean convex functions therefore show that the function $d_a$ restricted to 
$P \cup Q$ has Gateaux derivative $\frac{1}{\sqrt{5}}(1,2)$ at $x$.  The Euclidean normal cone to $P$ at $x$ is 
$\{0\} \times -\R_+$, so we deduce
\bmye \label{subdifferential}
\partial d_a^P(x) ~=~ \frac{1}{\sqrt{5}}(1,2) ~+~ (\{0\} \times -\R_+).
\emye
We now compare this with the subgradient provided by Theorem \ref{main}.  Let $y$ denote the point $(0,-1)$.  The geodesic from $x$ to $a$ consists of the two line segments $[x,y]$ in the square $Q$ and 
$[y,a]$ in the square $T$.  The common face of $P$ and $Q$ is $F = [0,1] \times \{0\}$, and in the square $Q$, the nearest point to $y$ in the face $F$ is the point $z = (0,0)$.  The unit vector in the direction $x-z$ is therefore $(1,0)$.  The angle $angle yxz$ has cosine $\frac{1}{\sqrt{5}}$, so the theorem asserts that the vector $\frac{1}{\sqrt{5}}(1,0)$ is a subgradient of $d_a^P$ at $x$.  This is confirmed by equation (\ref{subdifferential}).
}
\end{exa}

\section{Means, medians, and circumcenters}
In any metric space $(\X,d)$, we can consider a variety of interesting optimization problems associated with a given nonempty finite set $A \subset \X$, posed simply in terms of the distance functions $d_a$ to points 
$a \in A$.  Foremost among these are {\em weighted mean} problems, which entail minimizing functions $f \colon \X \to \R$ defined by
\bmye \label{weighted-mean}
f(x) ~=~ \sum_{a \in A} \lambda_a d_a^q(x) \qquad (x \in \X)
\emye
for some given weights $\lambda_a > 0$, for $a \in A$, and a given exponent $q \ge 1$.  The special case $q=1$ is known as the {\em median} problem.  Another well-known example is the {\em minimum enclosing ball} or {\em circumcenter} problem, which entails minimizing over $\X$ the function
\[
x ~\mapsto~ \max_{a \in A} d_a(x).
\]

When the space $\X$ is a finite CAT(0) cubical complex, we can solve weighted mean and minimum enclosing ball problems using Algorithm \ref{algorithm} and Theorem \ref{main}.  Consider first the weighted mean problem \cite[Example~6.3.4]{bacak-convex}.  The function (\ref{weighted-mean}) is convex, and a minimizer always exists, by \cite[Lemma 2.2.19]{bacak-convex}.  For $q>1$, the function $f$ is strictly convex, so the minimizer is unique.  In the case $q=1$, minimizers may not be unique.

To compute means and medians in general Hadamard spaces, the only methods previously analyzed appear to be the splitting proximal point algorithms described in \cite{bacak-medians}.  As we noted in the introduction, such methods are inevitably slow.  In the case of a finite CAT(0) cubical complex, we instead propose Algorithm \ref{algorithm} as a faster alternative.

To implement Algorithm \ref{algorithm}, we first fix attention on some cell $P$ in the complex $\X$.  Via an isometric embedding, we make the identification $P = [0,1]^n \subset \Rn$.  Using the notation of Theorem~\ref{main}, we seek to minimize the convex function
\[
f_P ~=~ \sum_{a \in A} \lambda_a (d_a^P)^q \colon P \to \R.
\]
For this purpose, we can use a standard Euclidean cutting plane method, of the kind described in \cite[Chapter~2]{bubeck}.  Such methods require a separation/subgradient oracle of the following form.  Consider any input point $x \in \Rn$.
\begin{itemize}
\item
If $x \not\in P$, then output a hyperplane separating $x$ from $P$.
\item
Otherwise, return the value $f_P(x)$ and a subgradient $v \in \partial f_P(x)$.
\end{itemize}
Theorem \ref{main} allows us to implement this oracle.  Separating any point from $P$ is trivial.  On the other hand, for any point $x$ in $P$, for each point $a \in A$, the theorem describes how to calculate a corresponding subgradient $v_a \in \partial d_a^P(x)$.  We deduce $qd_a^{q-1}(x)v_a \in  \partial \big((d_a^P)^q\big)(x)$, and then adding gives a subgradient:  
\[
\sum_a \lambda_a qd_a^{q-1}(x) v_a ~\in~ \partial f_P(x).
\]

\begin{exa}[A simple mean calculation, concluded]
{\rm
We return to the simple mean problem in Examples \ref{adriana} and \ref{simple2}.  On each of the three cells we solve the corresponding subproblem using a standard subgradient-based method, deploying the subgradients available through Theorem \ref{main}.  We compare six methods:  four cutting plane methods, the classical subgradient method, and the cyclic proximal point method from \cite{bacak-medians}.  The cutting plane methods consist of the classical ellipsoid algorithm and three standard bundle methods from \cite{preamb}:  a level bundle method, a proximal bundle method, and a dual level method.
The methods all converge to the optimal solution of each subproblem:  the true mean in the cell $P_1$, and the point $(0,0)$ in the other two cells.  The behavior of the objective value in the case of the cell $P_1$ is shown in Figure \ref{fig:subgradient}, plotted against the number of geodesics computed:  the behavior in the other cells is similar.  Even on this small low-dimensional example, the cyclic proximal point method is relatively slow.  The subgradient method, although clearly sublinear, converges reasonably quickly on this example, in part because the objective function happens to be smooth around the optimal solution.  The cutting plane methods are all faster, and the plot suggests the linear convergence we expect. 

\begin{figure}
    \centering
    \includegraphics[width=0.8\textwidth]{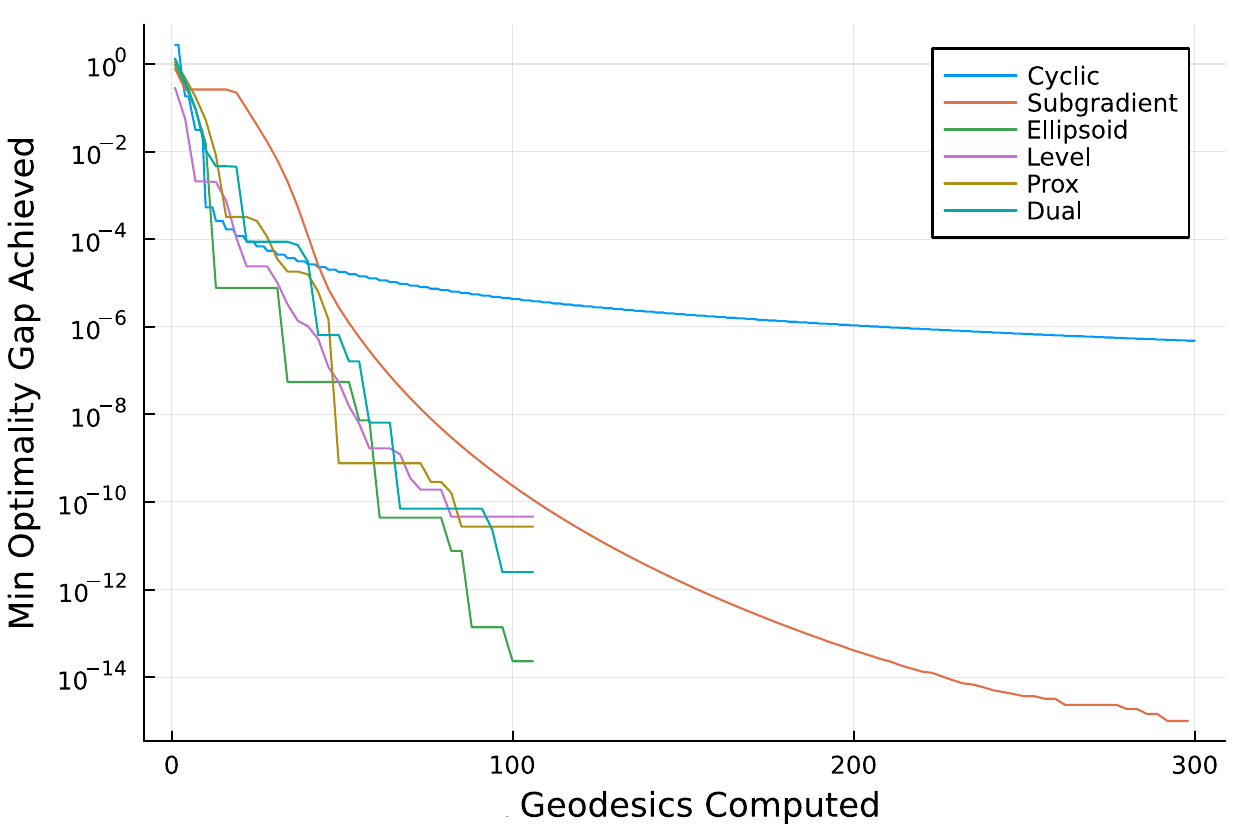}
    \caption{A simple mean calculation.  Five subgradient-based methods for minimizing, over a single cell, the sum of the squared distances to three points in a complex of three squares, compared with the cyclic proximal point method.  Each graph plots the smallest objective value found so far against the number of geodesics computed.
        }
    \label{fig:subgradient}
\end{figure}
}
\end{exa}

Turning to the minimum enclosing ball problem \cite[Example 2.2.18]{bacak-convex}, we can frame the problem as minimizing the strongly convex function $f \colon \X \to \R$ defined by
\[
f(x) ~=~ \max_{a \in A} d_a^2(x) \qquad (x \in \X).
\]
The unique solution is called the {\em circumcenter} of $A$.  Previous literature makes no apparent reference to algorithms for the minimum enclosing ball problem in general Hadamard spaces.  A subgradient-style method has been proposed recently in \cite{horoballs}, with the slow convergence rate typical of subgradient methods.   In the case of CAT(0) cubical complexes, we again propose Algorithm \ref{algorithm} as a faster alternative.

We follow the same approach, focusing on some cell $P$ and minimizing the strongly convex function $f_P \colon P \to \R$ defined by
\[
f_P(x)  ~=~ \max_{a \in A} (d_a^P)^2(x) \qquad (x \in \X).
\]
To compute a subgradient of $f_P$ at a point $x \in P$, we first select, from the set $A$, a point $a$ at maximum distance from $x$.  We then appeal to Theorem \ref{main} to calculate a corresponding subgradient $v_a \in \partial d_a^P(x)$.  Standard Euclidean subgradient calculus \cite{con_book} now gives us the desired subgradient:  $2d_a(x) v_a \in \partial f_P(x)$.

To apply Algorithm \ref{algorithm} to a mean, median, or circumcenter problem, having implemented the separation/subgradient oracle, we next choose a cutting plane method to solve the subproblems over cells $P$.  Most classically --- albeit not the best choice in theory or practice, as discussed in the introduction --- we could apply the ellipsoid algorithm, starting with the smallest Euclidean ball containing $P$.  That algorithm converges linearly:  by \cite[Theorem 2.4]{bubeck}, if the cell $P$ has dimension $n$, then after $t \ge n^2\log n$ calls to the oracle, we find an iterate in $P$ at which the value of $f$ exceeds $\min_P f$ by no more than 
\[
2\sqrt{n}(\max_P f) \exp \Big( -\frac{t}{2n^2} \Big).
\]
If $f$ is in fact {\em strongly} convex (as in the case of the mean or circumcenter problems), then the {\em iterates} also converge linearly, to the unique minimizer of $f$ on $P$.  

\section{The distance to an intersection of balls}
In addition to mean and circumcenter problems, another fundamental problem associated with a given nonempty finite subset $A$ of a metric space $(\X,d)$, and posed simply in terms of distance functions, is the {\em intersecting balls} problem.  Centered at each point $a \in A$, given a radius $\rho_a > 0$, we consider the closed ball $B_{\rho_a}(a)$.  We seek a point in their intersection: 
\[
\F = 
\bigcap_{a \in A} B_{\rho_a}(a).
\]

In any Hadamard space $\X$, we could approach the intersecting balls problem via the method of cyclic projections \cite{sims-cyclic,genaro-cyclic}:  the projection of any point $x$ onto a ball centered at $a$  is easy to compute from the geodesic $[x,a]$, since the projection must lie on the geodesic.  However, in a general Hadamard space, the rate at which the method converges to a point in the intersection is not clear, and when the intersection is empty, the behavior of the iterates is not well understood \cite{lytchak-cyclic}.  

When the space $\X$ is a finite CAT(0) cubical complex, we take a different approach to the intersecting balls problem, again using Algorithm \ref{algorithm} and Theorem~\ref{main}.  Consider the related problem of minimizing over $\X$ the objective function
\[
\sum_{a \in A}  \mbox{dist}_{B_{\rho_a}(a)}.
\]
Since this function is coercive, it has a minimizer, by \cite[Lemma~2.2.19]{bacak-convex}.  Consider any such minimizer $\bar x$.  If $\bar x \in \F$, then we have solved our problem, and on the other hand, if $\bar x \not\in \F$, then the set $\F$ must be empty.  We therefore seek to minimize the function $f \colon \X \to \R$ defined by
\[
f(x) ~=~ \sum_{a \in A}  \max\{d_a(x) - \rho_a, 0\} \qquad (x \in \X).
\]

We follow the same strategy as for the mean problem in the previous section, focusing on individual cells $P$ in the complex $\X$, making the identification $P = [0,1]^n$ and seeking to minimize the convex function $f_P \colon P \to \R$ defined by
\[
f_P(x) ~=~  \sum_{a \in A}  \max\{d_a^P(x) - \rho_a, 0\} \qquad (x \in P).
\]
As before, to implement a cutting plane algorithm, we need to find a subgradient of $f_P$ at any input point $x \in P$.  To that end, for each point $a \in A$, we compute the distances $d_a(x)$, and calculate a vector $v_a$ as follows.  If $d_a(x) \le \rho_a$, then we define $v_a = 0$.  On the other hand, if 
$d_a(x) > \rho_a$, then we appeal to Theorem \ref{main} to calculate a subgradient $v_a \in \partial d_a^P(x)$.  In that way, standard Euclidean subgradient calculus \cite{con_book} guarantees that $v_a$ is a subgradient of the convex function given by
\[
x ~\mapsto~  \max\{d_a^P(x) - \rho_a, 0\} \qquad (x \in P),
\]
from which we deduce our desired subgradient:
\[
\sum_{a \in A} v_a ~\in~ \partial f_P(x).
\]
We can now apply a cutting plane algorithm.

Refining the intersecting balls problem, we might seek to compute the distance from a given point $b \in \X$ to the intersection of balls $\F$.  We can approach this problem by a bisection search strategy.  We first run the previous algorithm, either finding a point $\hat x \in \F$ or detecting that none exists, in which case we terminate.  If we find $\hat x \in \F$, we begin our bisection search, always maintaining an interval
$[\underline{\rho},\overline{\rho}]$ containing the distance from $b$ to $\F$.  Initially we set 
$\underline{\rho} = 0$ and $\overline{\rho} = d(b,\hat x)$.  At each iteration, we consider the midpoint 
$\rho$ of the current interval, and apply the previous algorithm to the intersecting balls problem 
$\F \cap B_{\rho}(b)$.  If the intersection is empty, we update $\underline{\rho} = \rho$;  otherwise we update $\overline{\rho} = \rho$.  We then repeat.

\small
\parsep 0pt

\def\cprime{$'$} \def\cprime{$'$}

\appendix
\normalsize
\section{Appendix:  core cells} \label{appendix}
An interesting source of illustrative computational examples arises from the class of cubical complexes $\X$ consisting of a collection of cells all sharing one common face.  As observed in \cite{sticky}, although this face is lower-dimensional than $\X$ it often contains the mean of a set of points distributed throughout $\X$.  

\begin{exa}[An open book \cite{sticky}]
{\rm
Consider three copies of the square $[0,1]^2$ joined along the shared edge $\{0\} \times [0,1]$:  the ``spine'' of the open book.  Denote the three copies $P^1,P^2,P^3$, and consider any points $a^i \in \inT P^i$ satisfying
\[
a^i < \sum_{j \ne i} a^j \qquad \mbox{for}~ i=1,2,3.
\]
(An example is $a^i = (\frac{1}{2},\frac{1}{2})$ for each $i$.)
A quick calculation shows that the mean is the following point on the spine:
\[
\Big(0, \frac{1}{3}(a^1_2+ a^2_2 + a^3_2)\Big).  
\]
The spine is thus sticky:  the mean remains there for all small perturbations of the three points.
}
\end{exa}

\begin{defn}
{\rm
In a cubical complex, if all distinct pairs of maximal cells have the same nonempty intersection $F$, then we call $F$ the {\em core}.
}
\end{defn}

\noindent
The core, if it exists, is a common face of all the maximal cells.  For example, the 3-spider has the core $\{0\}$.  The maximal cells may not all have the same dimension.  For example, the complex consisting of the maximal cells in $\R^3$
\[
[0,1] \times [0,1] \times [0,1] \quad \mbox{and} \quad [-1,0] \times \{0\} \times [0,1]
\quad \mbox{and} \quad \{0\} \times [-1,0]  \times [0,1]
\]
has core $\{0\} \times \{0\} \times [0,1]$.

\begin{prop}
Any cubical complex with a core is {\rm CAT(0)}.
\end{prop}

\pf
Consider a cubical complex $\X$ with a core $F$. The cells of $\X$ consist just of all the faces of maximal cells.  Since $\X$ is contractible to any point in $F$, it is simply connected.

We next observe that Gromov's link condition holds.  To see this, given any vertex, consider any finite set $E$ of $k$ incident edges, each pair of which lie in a common $2$-dimensional cell of $\X$.  For each maximal cell $P$, let $E_P$ denote those edges in $E$ that lie in $P$ but not in the common face $F$.  Consider any two maximal cells $P,Q$.  Suppose that there exist edges $p \in E_P \setminus E_Q$ and $q \in E_Q \setminus E_P$.  Any $2$-dimensional cell containing both $p$ and $q$ must be a face of some maximal cell $R$.  The edge $p$ lies in both $R$ and $P$, but by assumption not in $F$, so in fact $R=P$.  But the same argument shows $R=Q$, giving a contradiction.  We deduce that either $E_P \subset E_Q$ or $E_Q \subset E_P$, proving that the collection of sets $E_P$ for maximal cells $P$ is totally ordered by inclusion.  The collection therefore has maximal element, $E_P$ for some maximal cell $P$.  Every edge in $E$ lies in this cell $P$, and hence in some $k$-dimensional face of $P$, as required.
\finpf

We next compute geodesics.

\begin{prop}
Consider points $x$ and $y$ in a cubical complex $\X$ with a core $F$.  If some cell contains both points, then the geodesic $[x,y]$ is the line segment joining them in that cell.  On the other hand, if no cell contains both,  then there are maximal cells $P$ and $Q$ such that $x \in P \setminus Q$ and 
$y \in Q \setminus P$.  In $P$, consider the distance $|x-x_F|$ between $x$ and $x_F$, its orthogonal projection onto $F$.  Analogously, in $Q$, consider the distance $|y-y_F|$ between $y$ and $y_F$, its orthogonal projection onto $F$.  Consider the following convex combination in $F$:
\[
z ~=~ \frac{1}{|x-x_F| + |y-y_F|} \big(|y-y_F|x_F ~+~ |x-x_F|y_F \big).
\]
Then the geodesic $[x,y]$ is the union of the two geodesics $[x,z]$ and $[z,y]$.
\end{prop}

\pf
Note that the geodesic $[x,y]$ must pass through the face $F$, so it is the union of two geodesics $[x,w]$ and $[w,y]$ for some point $w \in F$.  The length of this path is
\[
\sqrt{|w-x_F|^2 + |x-x_F|^2} ~+~ \sqrt{|w-y_F|^2 + |y-y_F|^2}.
\]
As $w$ varies over $F$, this function is convex, and viewed on any ambient Euclidean space in which $F$ lives, its gradient is
\bmye \label{gradient}
\frac{1}{\sqrt{|w-x_F|^2 + |x-x_F|^2}}(w-x_F)
~+~
\frac{1}{\sqrt{|w-y_F|^2 + |y-y_F|^2}}(w-y_F).
\emye
Notice 
\[
z-x_F ~=~ 
\frac{|x-x_F|}{|x-x_F| + |y-y_F|}(y_F - x_F),
\]
so
\[
\sqrt{|z-x_F|^2 + |x-x_F|^2}
\: = \:
\frac{|x-x_F|}{|x-x_F| + |y-y_F|} \sqrt{|y_F - x_F|^2 + (|x-x_F| + |y-y_F|)^2}
\]
and hence the first term in the expression (\ref{gradient}), when $w=z$, becomes
\[
\frac{1}{\sqrt{|y_F - x_F|^2 + (|x-x_F| + |y-y_F|)^2}}
(y_F - x_F).
\]
By symmetry, the second term is identical, after multiplying by $-1$.  Hence the gradient (\ref{gradient})
is zero when $w=z$, so $z$ indeed minimizes the length.
\finpf


\end{document}